\documentstyle{amsppt}
\tolerance 3000 \pagewidth{5.5in} \vsize7.0in
\magnification=\magstep1 \widestnumber \key{AAAAAAAAAAA}
\NoRunningHeads \topmatter
\author S. Hofmann and A. Iosevich
\endauthor
\title Falconer conjecture in the plane for random metrics
\endtitle
\abstract The Falconer conjecture says that if a compact planar
set has Hausdorff dimension $>1$, then the Euclidean distance set
$\Delta(E)=\{|x-y|: x,y \in E\}$ has positive Lebesgue measure. In
this paper we prove, under the same assumptions, that for almost
every ellipse $K$, $\Delta_K(E)=\{{||x-y||}_K: x,y \in E\}$ has
positive Lebesgue measure, where ${||\cdot||}_K$ is the norm
induced by an ellipse $K$. Equivalently, we prove that if a
compact planar set has Hausdorff dimension $>1$, then $\Delta(TE)$
has positive Lebesgue measure for almost every transformations $T$
with bounded positive eigenvalues. We also use this result to
deduce a version of the Erdos Distance Conjecture in the plane.
\endabstract
\date December 19, 2002
\enddate
\address S. Hofmann, University of Missouri
\ email: hofmann \@ math.missouri.edu
\endaddress
\address A. Iosevich, University of Missouri
\ email: iosevich \@ math.missouri.edu
\endaddress
\thanks The research of the authors is partially supported by NSF
grants \endthanks
\endtopmatter
\document

\head Introduction \endhead

\vskip.125in

Let $E \subset {[0,1]}^d$, $d \ge 2$. The celebrated Falconer
conjecture says that if the Hausdorff dimension of $E$ exceeds
$\frac{d}{2}$, then the distance set $\Delta(E)=\{|x-y|: x,y \in
E\}$ has positive Lebesgue measure.

The initial result in this direction was proved by Falconer
(\cite{Falconer86}) who showed that $\Delta(E)$ has positive
Lebesgue measure if the Hausdorff dimension of $E$ exceeds
$\frac{d+1}{2}$. This result was later improved in all dimensions
by Bourgain (\cite{Bourgain94}). The best known result in the
plane is due to Tom Wolff who proved that $\Delta(E)$ has positive
Lebesgue measure provided that the Hausdorff dimension of $E$ is
greater than $\frac{4}{3}$.

The purpose of this paper is to prove that the conclusion of the
Falconer conjecture holds for almost every linear perturbation of
the Euclidean metric. More precisely, let
$$\Delta_K(E)=\{{||x-y||}_K: x,y \in E\}, \tag0.1$$ where $K$ is a
symmetric bounded convex set in ${\Bbb R}^2$ and ${||\cdot||}_K$
is the distance induced by $K$. Our main result is the following.

\proclaim{Theorem 0.1} Let $E \subset {[0,1]}^2$ be a set of
Hausdorff dimension greater than $1$. Let $K_{a,\phi}$ denote the
ellipse with eccentricities $a_1, a_2$, $1 \leq a_j \leq 2$,
rotated by the angle $\phi$, and let $\Delta_{a,\phi}(E)$ denote
the corresponding distance set. Then ${\Bbb
E}(m(\Delta_{a,\phi}(E)))>0$, where the expectation is taken with
respect to the uniform distribution on $[1,2] \times [1,2] \times
[0,\pi] $, and $m$ denotes the one-dimensional Lebesgue measure.
\endproclaim

The sharpness of Theorem 0.1 is demonstrated by a modification of
a construction due to Falconer (\cite{Falconer86}). Let $0<s \leq
2$. Let $q_1, q_2, \dots, q_i \dots$ be a sequence of positive
integers such that $q_{i+1} \ge q_i^i$. Let $E_i=\{x \in {\Bbb
R}^2: 0 \leq x_j \leq 1, |x_j-p_j/q_i| \le q_i^{-\frac{2}{s}} \
\text{for some integers} \ p_j, \ j=1,2\}$. It is not hard to see
(see e.g. \cite{Wolff02}) that the Hausdorff dimension of
$E=\cap_{i=1}^{\infty} E_i$ is $s$. Also, $\Delta(E) \subset
\cap_{i=1}^{\infty} \Delta(E_i)$.

Let $P_i=\{p=(p_1,p_2): 0 \leq p_j \leq q_i\}$. Let
$\Delta^{(0,0)}_{a, \phi}(P_i)=\{{||p||}_{a,\phi}: p \in P_i\}$.
It is immediate that $\# \Delta^{(0,0)}_{a, \phi}(P_i) \leq
{(q_i+1)}^2$. By translation invariance it follows that $\#
\Delta_{a,\phi}(P_i) \leq {(q+1)}^2$. We conclude that
$\Delta_{a,\phi}(E_i)$ is contained in at most ${(q_i+1)}^2$
intervals of length $\approx q_i^{-\frac{2}{s}}$. It follows that
the Hausdorff dimension of $\Delta(E)$ is $\leq s$. Thus if $s<1$,
$\Delta_{a,\phi}(E)$ has Lebesgue measure $0$ for every $a, \phi$.

Observe that in the above example it is not necessary for $p_i$s
to be integers. It is quite sufficient for $P_i$ to be
sufficiently dense, and separated in the sense that there exists
$c>0$ such that $|p-p'| \ge c$, $p,p' \in P_i$, $p \not=p'$. This
observation allows us to use Theorem 0.1 to deduce a version of
the Erdos Distance Conjecture. See e.g. \cite{PaAg95} for a
thorough description of the Erdos Distance Problem and related
concepts.

The classical Erdos Distance Problem is to obtain a lower bound
for $\Delta(S)$, where $S$ is a finite subset of ${\Bbb R}^2$.
Erdos Distance Conjecture says that for any $\epsilon>0$ there
exists $C_{\epsilon}>0$ such that $\# \Delta(S) \ge
C_{\epsilon}{(\# S)}^{1-\epsilon}$. A slightly weaker version of
the Erdos Distance can be stated as follows.

\proclaim{Asymptotic version of the Erdos Distance Conjecure} Let
$A$ be a separated subset of ${\Bbb R}^2$. Suppose that $A$ is
actually a Delone set, which means that there exists a universal
constant $C>0$ such that the intersection of $A$ with any cube of
side-length $R$ contains $\ge CR^2$ elements of $A$. Then for any
$\epsilon>0$ there exists a positive constant $C'_{\epsilon}>0$
such that
$$ \# \Delta(A \cap {[-R,R]}^2) \ge
C'_{\epsilon} R^{2-\epsilon}. \tag0.1$$
\endproclaim

Using the above counter-example used to establish sharpness of
Theorem 0.1 we can prove the following random variant of the
asymptotic version of the Erdos Distance Conjecture. See
\cite{IoLa2003} for a systematic application of this mechanism to
non-Euclidean distances in ${\Bbb R}^d$.

\proclaim{Corollary 0.2} Let $A$ be as in the statement of the
asymptotic version of the Erdos Distance Conjecture. Then for any
$\epsilon>0$, there exists $C'_{\epsilon}>0$ such that
$$ \# \Delta_{a,\phi}(A \cap {[-R,R]}^2) \ge
C'_{\epsilon} R^{2-\epsilon}. \tag0.2$$ for almost every $(a,
\phi) \in {[1,2]}^2 \times [0, \pi]$.
\endproclaim

To prove Corollary 0.2, let $E_i, E,$ and $\{q_i\}$ be defined as
in the counter-example above with $P_i=A \cap {[0,q_i]}^2$.
Suppose that $\# \Delta_{a,\phi}(P_i) \leq Cq_i^{2-\epsilon}$ for
some $\epsilon>0$ for a sequence of $i$s going to infinity. Then
$\Delta_{a,\phi}(E_i)$ can be covered by $\leq Cq_i^{2-\epsilon}$
intervals of length $\approx q^{-\frac{2}{s}}$. It follows that
the Hausdorff dimension of $\Delta_{a, \phi}(E)$ is $\leq
s-\frac{s\epsilon}{2}$. Now let $s=1+\delta$, $\delta>0$. We
conclude that the Hausdorff dimension of $\Delta_{a, \phi}(E)$ is
$\leq 1+\delta-\frac{(1+\delta)\epsilon}{2}<1$ if $\delta$ is
sufficiently small. This is a contradiction because Theorem 0.1
implies that $\Delta_{a, \phi}(E)$ has positive Lebesgue measure
for almost every $(a, \phi)$.

\remark{Remark 1} Another way of stating Theorem 0.1 is to say
that if $E \subset {[0,1]}^2$ be a set of Hausdorff dimension
greater than $1$, then for almost every linear transformation $T$
of the form rotation followed by an anisotropic dilation,
$\Delta(TE)$ has positive Lebesgue measure. \endremark

\remark{Remark 2} The proof of Theorem 0.1 below will show that
the conclusion of Theorem 0.1 still holds if in the definition of
$\Delta_{a, \phi}$, the Euclidean circle is replaced by any smooth
curve with everywhere non-vanishing curvature. In other words,
Theorem 0.1 does not just hold for perturbations of the Euclidean
metric, but also for perturbations of any metric whose unit circle
is smooth and has non-vanishing curvature. \endremark

\remark{Remark 3} In principle, an appropriate variant of Theorem
0.1 should hold in the context of two-dimensional Riemannian
manifolds. We shall address this issue in a subsequent paper.
\endremark

\remark{Remark 4} It is worth noting that large classes of
two-dimensional sets of Hausdorff dimension $\alpha>1$ for which
the Falconer conjecture holds can be constructed using more
complicated probabilistic schemes. For example, let $A$ be a
compact subset of the real line of Hausdorff dimension
$\frac{1}{2}<\frac{\alpha}{2}<1$, and $W(t)$ an almost surely
continuous version of the real Wiener process. A theorem due to
J.P. Kahane (see e.g. \cite{Kahane68}) says that $E=W(A)$ is
almost surely a Salem set of dimension $\alpha>1$, which means
that $E$ is equipped with a Borel measure $\mu$ such that
$|\widehat{\mu}(\xi)| \leq C{(1+|\xi|)}^{-\frac{\alpha}{2}}$.
Using Theorem 1.1 below one easily deduces that $\Delta(E)$ has
positive Lebesgue measure, so Falconer conjecture holds for this
class of fractal sets. \endremark

\remark{Remark 5} In particular, the proof of Corollary 0.2 shows
that Falconer Distance Conjecture implies the asymptotic Erdos
Distance conjecture. It would be nice to prove that the Falconer
conjecture in fact implies the standard Erdos Distance Conjecture.
This amounts to eliminating the well-distributivity assumption on
$A$ in the statement of Corollary 0.2 and replacing $C'_{\epsilon}
R^{2-\epsilon}$ on the right hand side of $(0.2)$ by
$C'_{\epsilon} {(\# A \cap {[-R,R]}^2)}^{1-\epsilon}$.
\endremark

\vskip.125in

\head Method of proof of Theorem 0.1 \endhead

\vskip.125in

We use a modification of the following result due to Mattila
(\cite{Mattila87}).

\proclaim{Theorem 1.1} Let $E \subset {[0,1]}^2$ with a Borel
measure $\mu$. Suppose that
$$ \int_1^{\infty} {\left( \int_{S^1} {|\widehat{\mu}(t\omega)|}^2
d\omega \right)}^2 t dt<\infty. \tag1.1$$

Then $\Delta(E)$ has positive Lebesgue measure. (Here and
throughout $\Delta(E)=\Delta_K(E)$ with $K$ a unit disk).
\endproclaim

In fact, the argument used to prove Theorem 0.1 combined with a
standard stationary phase argument (see e.g. Theorem 0.4 below)
yields the following slightly more general result:
\proclaim{Theorem 1.2} Let $E \subset {[0,1]}^2$ with Borel
measure $\mu$. Let $K$ be a bounded convex set such that $\partial
K$ is smooth and has everywhere non-vanishing curvature. Suppose
that
$$ \int_1^{\infty} {\left( \int_{\partial K}
{|\widehat{\mu}(t\omega_K)|}^2 d\omega_K \right)}^2 t dt<\infty,
\tag1.2$$ where $d\omega_K$ denotes the Lebesgue measure on
$\partial K$, the boundary of $K$.

Let $K^{*}=\{\xi: \sup_{x \in K} x \cdot \xi \leq 1\}$, the convex
set dual to $K$. Then $\Delta_{K^{*}}(E)$ has positive Lebesgue
measure. \endproclaim

We shall give a proof of Theorem 1.2 at the end of this paper for
the sake of completeness.

In view of Theorem 1.1 and Theorem 1.2, Theorem 0.1 follows from
the following estimate

\proclaim{Theorem 1.3} Let $E \subset {[0,1]}^2$ with Borel
measure $\mu$. Suppose that the Hausdorff measure of $E$ is
greater than $1$. Then
$$ \int_0^{\pi} \int \int_1^{\infty} {\left( \int_{S^1}
{|\widehat{\mu}(t\omega_{a, \phi})|}^2 d\omega \right)}^2 t dt
\psi(a)da d\phi <\infty, \tag1.3$$ where $\omega_{a, \phi}$ is the
standard parameterization of the ellipse with eccentricities $a_1,
a_2$, rotated by $\phi$, and $\psi$ is a smooth cutoff function
identically equal to $1$ in ${[1,2]}^2$ and vanishing outside
${[1/2,4]}^2$.
\endproclaim

The result due to Wolff mentioned above was proved by showing that
under the assumptions of Theorem 0.1,
$$ \int_{S^1} {|\widehat{\mu}(t \omega)|}^2 d\omega \lesssim
t^{-\frac{\alpha}{2}}, \tag1.4$$ and an example due to Sjolin
(\cite{Sjolin93}) shows that this estimate cannot, in general, be
improved. This means that the proof of Theorem 0.3 must heavily
rely on averaging in $t$, $a$, and $\phi$.

Throughout the paper we shall make use of the following version of
the method of stationary phase. See e.g. \cite{Sogge93}, Theorem
1.2.1.

\proclaim{Theorem 1.4} Let $S$ be a convex smooth hyper-surface in
${\Bbb R}^d$ with everywhere non-vanishing Gaussian curvature and
$d\mu$ a $C_0^{\infty}$ measure on $S$. Then
$$ |\widehat{d\mu}(\xi)| \lesssim {|\xi|}^{-\frac{d-1}{2}}.
\tag1.5$$

Moreover, suppose that $\Gamma \in {\Bbb R}^d \ \backslash \ (0,
\dots, 0)$ is the cone consisting of all vectors $\xi$ normal to
$S$ at some point $x$ in a fixed relatively compact neighborhood
of support of $d\mu$. Then
$$ \left|{\left( \frac{\partial}{\partial \xi}\right)}^{\alpha}
\widehat{d\mu}(\xi) \right|=O({(1+|\xi|)}^{-N}) \ \forall \ N, \
\text{if} \ \xi \notin \Gamma, \tag1.6$$ and
$$ \widehat{d\mu}(\xi)=\sum_{j=1}^2 e^{-2 \pi i x_j \cdot \xi}
a_j(\xi), \ \text{if} \ \xi \in \Gamma, \tag1.7$$ where the finite
sum is taken over the points $x_j \in {\Cal N}$ having $\xi$ as a
normal and
$$ \left|{\left( \frac{\partial}{\partial
\xi}\right)}^{\alpha}a_j(\xi) \right| \leq
C_{\alpha}{(1+|\xi|)}^{-\frac{d-1}{2}-|\alpha|}. \tag1.8$$
\endproclaim

{\bf Notation:} Throughout this paper, $a \lesssim b$ means that
there exists a positive constant $C$ such that $a \leq Cb$.
Similarly, $a \lessapprox b$, with respect to a parameter $s$,
means that given $\epsilon>0$ there exists $C_{\epsilon}>0$ such
that $a \leq C_{\epsilon}s^{\epsilon} b$.

\vskip.125in

\head Proof of Theorem 1.3 \endhead

\vskip.125in

Let $\beta \in C_0^{\infty}$ be the usual Littlewood-Paley cutoff,
i.e $\beta$ is supported in $[1/2,4]$, $\beta \equiv 1$ in
$[1,2]$, and $\sum_j \beta(2^{-n} \cdot) \equiv 1$. Let
$\omega_{a, -\phi}=\rho_{\phi}^{-1} \delta_a \omega$, where
$\rho_{\phi}$ denotes the rotation by the angle $\phi$, and
$\delta_a(x)=(a_1x_1,a_2x_2)$. Let $\mu$ be a probability measure
on $E$ such that $\mu(\{x \in E: |x-y| \leq r\}) \leq
Cr^{\alpha}$, $r>0$, where $\alpha$ is the Hausdorff dimension of
$E$. For the existence of such a measure, see, for example,
Proposition 8.2 in \cite{Wolff02}. Define
$$ I_n=\int \int {\left( \int_{S^1}
{|\widehat{\mu}(t\omega_{a,-\phi})|}^2 d\omega \right)}^2 t
\beta(2^{-n}t) dt \psi(a)da d\phi$$
$$=\int \int \int \int {|\widehat{\mu}(t\omega_{a,-\phi})|}^2
{|\widehat{\mu}(t\omega'_{a,-\phi})|}^2 d\sigma(\omega)
d\sigma(\omega') t \beta(2^{-n}t) dt \psi(a)da d\phi$$
$$=\int \int \int \int \int e^{2 \pi i ((x-y) \cdot
t \omega'_{a,-\phi}+(x'-y') \cdot t \omega'_{a,-\phi})} d\mu^{*}
d\sigma(\omega) d\sigma(\omega')t \beta(2^{-n}t) dt \psi(a)da
d\phi, \tag2.1$$ where
$$ d\mu^{*}=d\mu(x)d\mu(y)d\mu(x')d\mu(y'), \tag2.2$$ and
$d\sigma$ denotes a $C_0^{\infty}$ measure on the sphere. Using a
partition of unity we see that it is enough to consider this
situation.

Integrating in $\omega$ and $\omega'$ first, we get
$$ \int \int \int \int \widehat{\sigma}(t{(x-y)}_{a,\phi})
\widehat{\sigma}(t{(x'-y')}_{a,\phi})d\mu^{*}\beta(2^{-n}t) dt
\psi(a)da d\phi, \tag2.3$$ where $\sigma$ is $C_0^{\infty}$
measure on $S^1$ as above, and
$x_a=(a_1(x_1\cos(\phi)-x_2\sin(\phi)),
a_2(x_1\sin(\phi)+x_2\cos(\phi)))$.

\subhead Case 1: $2^n|x-y| \lesssim 1$ and $2^n|x'-y'| \lesssim 1$
\endsubhead Then for any $\epsilon>0$, $(1.3)$ is bounded by
$$ \int {|x-y|}^{-\alpha+\epsilon} {|x'-y'|}^{-\alpha+\epsilon}
t^{-2 \alpha+2\epsilon} d\mu(x)d\mu(y)d\mu(x') d\mu(y') tdt<\infty
\tag2.4$$ as desired if $\epsilon$ is sufficiently small.

\subhead Case 2: $2^n|x-y|>>1$ and $2^n|x'-y'|>>1$ \endsubhead
Observe that the symbol of order $0$ resulting from pulling
${(t|x-y|)}^{-\frac{1}{2}}$ from the symbol $a_j$ given by Theorem
1.4 can be incorporated into the smooth cut-off $\beta$ without
effecting the size or the support of $\beta$ or its derivatives.
We shall suppress (harmless) dependence of $\beta$ on
$x,y,x',y',a, \phi$ in what follows. Using this observation and
$(1.7)$ above, we see that $(2.3)$ can be written as a sum of
terms of the form
$$ \int \int \int \int e^{2 \pi i
t({|x-y|}_{a,\phi}-{|x'-y'|}_{a,\phi})} t^{-1} {|x-y|}_{a,
\phi}^{-\frac{1}{2}} \ {|x'-y'|}_{a, \phi}^{-\frac{1}{2}} \ t
d\mu^{*}\beta(2^{-n}t) dt \psi(a)da d\phi. \tag2.5$$

We must also consider the term where the phase function is $2 \pi
i t({|x-y|}_{a,\phi}+{|x'-y'|}_{a,\phi})$, but this case is very
easy. Let $\eta$ be a small parameter to be determined later. We
have
$$ \int \int \int \int e^{2 \pi i
t({|x-y|}_{a,\phi}+{|x'-y'|}_{a,\phi})} t^{-1} {|x-y|}_{a,
\phi}^{-\frac{1}{2}} \ {|x'-y'|}_{a, \phi}^{-\frac{1}{2}} \ t
d\mu^{*}\beta(2^{-n}t) dt \psi(a)da d\phi$$
$$ =2^n \int_{\{(x,y,x',y'): ||x-y|+|x'-y'|| \lesssim
2^{-n(1-\eta)}\}} {|x-y|}^{-\frac{1}{2}} {|x'-y'|}^{-\frac{1}{2}}
d\mu^{*}+O(2^n 2^{-n \eta N})$$
$$ \lesssim 2^n \int {|x-y|}^{-\alpha+\epsilon}
2^{-n(1-\eta)(\alpha-\frac{1}{2}-\epsilon)}
{|x'-y'|}^{-\alpha+\epsilon}
2^{-n(1-\eta)(\alpha-\frac{1}{2}-\epsilon)} d\mu^{*}$$ $$ \lesssim
2^{2n(1-\eta)} 2^{-2n(1-\eta) \alpha} 2^{2n(1-\eta) \epsilon} \tag
2.6$$ which sums if $\epsilon<\alpha-1$ and $\eta<1$. The second
line of $(2.6)$ follows from the first using the fact, which
follows easily by integration by parts, that the Fourier transform
of a smooth compactly supported function decays rapidly at
infinity.

We now turn our attention to $(2.5)$. Integrating in $t$ first we
get
$$ 2^n \int \int \int
\widehat{\beta}(2^n({|x-y|}_{a,\phi}-{|x'-y'|}_{a,\phi}))
{|x-y|}_{a, \phi}^{-\frac{1}{2}} \ {|x'-y'|}_{a,
\phi}^{-\frac{1}{2}}\psi(a)da d\phi d\mu^{*}. \tag2.7$$

Localizing to the sets where $2^{-k} \leq {|x-y|}_{a,\phi} \leq
2^{-k+1}$, $2^{-k'} \leq {|x'-y'|}_{a,\phi} \leq 2^{-k'+1}$, we
obtain
$$ I_{n,k,k'} \approx 2^n 2^{\frac{k}{2}} 2^{\frac{k'}{2}}
\int \int \int
|\widehat{\beta}(2^n({|x-y|}_{a,\phi}-{|x'-y'|}_{a,\phi}))|
\psi(a)da d\phi d\mu^{*} \tag2.8$$

Let $x-y=|x-y|e^{iA}$, $x'-y'=|x'-y'|e^{iB}$. We now decompose
$x-y$ and $x'-y'$ into sectors of aperture $\delta$ to be
determined. Let $S^{j,k}_{\delta}$ denote the "rectangle" formed
by intersection the annulus $\{z: 2^{-k} \leq |z| \leq 2^{-k+1}\}$
and the angular sector $\{z: j\delta \leq A \leq (j+1)\delta\}$.
Define $S_{\delta}^{j',k'}$ analogously.

Let
$$ I_{n,k,k'}^{j, j', \delta}=2^n 2^{\frac{k}{2}} 2^{\frac{k'}{2}}
\int_{S_{\delta}^{j,k} \bigotimes S_{\delta}^{j',k'}} \int \int
|\widehat{\beta}(2^n({|x-y|}_{a,\phi}-{|x'-y'|}_{a,\phi}))|
\psi(a)da d\phi d\mu^{*}$$
$$ =2^{n+\frac{k}{2}+\frac{k'}{2}}
\int_{\{x-y \in S_{\delta}^{j,k}; x'-y' \in S_{\delta}^{j',k'};
|{|x-y|}_{a,\phi}-{|x'-y'|}_{a,\phi}| \lesssim 2^{-n(1-\eta)}\}}
\int \int \psi(a)da d\phi d\mu^{*}$$ $$+O(2^n 2^{\frac{k}{2}}
2^{\frac{k'}{2}}2^{-n \eta N}), \tag2.9$$ where $\eta>0$ is a
small parameter to be chosen later.

Observe that for $j, j'$ fixed, we have that $x-y$ is in a $\delta
2^{-k}$ by $2^{-n(1-\eta)}$ rectangle and $x'-y'$ is in a $\delta
2^{-k'}$ by $2^{-n(1-\eta)}$ rectangle. Also observe that if both
$k,k' \gtrsim n(1-\eta)$, then we have a simple estimate analogous
to the one in Case 1 if $\eta$ is chosen to be sufficiently small
and $\delta \approx 1$. Otherwise, if at least one of
$k,k'<<n(1-\eta)$, then they both are, and, moreover, $k \approx
k'$. Therefore, in what follows we may assume that we are in the
latter situation, so that the double index appearing above may now
be replaced by the single index $k$.

Now, multiplying both sides by
$|{|x-y|}_{a,\phi}+{|x'-y'|}_{a,\phi}|$ and computing the area of
the resulting set, we see that
$$ |\{(a_1,a_2): |{|x-y|}_{a,\phi}-{|x'-y'|}_{a,\phi}|
\lesssim 2^{-n(1-\eta)} \}| \lesssim 2^{-(n(1-\eta)-k)}
{[|A-B||A+\phi|]}^{-1} \tag2.10$$ where, without loss of
generality, $A \geq B$, so that $j-j' \ge 1$ unless $j=j'$. We
also take $A,B, \phi$ to be small and positive. (The other cases
follow by the same argument). It follows that if $j \not=j'$,
$$ I_{n,k}^{j,j',\delta} \lesssim 2^{n+k}
\int_S \int \max \left\{1, \
\frac{2^{-(n(1-\eta)-k)}}{|j-j'|\delta|\delta j+\phi|} \right\}
d\phi d\mu^{*}. \tag2.11$$ where integration is over the set $S$
where $x-y \in S_{\delta}^{j,k}, x'-y' \in S_{\delta}^{j',k'}$
(recall that $k \approx k'$). Let $\delta \approx
2^{-n(1-\eta)+k}$. Then for $j'$ and $x-y$ fixed, we have that
$x'-y'$ is located in a ball of radius $\approx 2^{-n(1-\eta)}$.
We have
$$ I_{n,k}^{j,j'} \lesssim 2^{n+k} \int_S \int \max \left\{1,
\ \frac{1}{|j-j'||\delta j+\phi|} \right\} d\phi d\mu^{*}$$
$$ \lessapprox 2^{n+k} \frac{1}{|j-j'|} \int_S
d\mu(x)d\mu(y)d\mu(x')d\mu(y'). \tag2.12$$

We must estimate
$$ \sum_{k,n} 2^{n+k} \sum_{j,j'} \frac{1}{|j-j'|} \int_S
d\mu(x)d\mu(y)d\mu(x')d\mu(y'). \tag2.13$$

Let $l=j-j'$. We get
$$ \sum_{k,n} 2^{n+k} \sum_{l} \frac{1}{l} \sum_{j'} \int_S
d\mu(x)d\mu(y)d\mu(x')d\mu(y'). \tag2.14$$

For a fixed $x-y$ and a sector given indexed by $j$, $x'-y'$ is
contained in a ball of radius $C2^{-n(1-\eta)}$ since $\delta
2^{-k} \approx 2^{-n(1-\eta)}$. Fixing $y'$ and integrating in
$x'$ we get $2^{-n(1-\eta) \alpha}$ since $\mu$ is
$\alpha$-dimensional. Taking the union over all the sectors
indexed by $j'$, we have $x-y$ in the annulus of width $2^{-k}$.
Fixing $y$ and integrating in $x$, we pick up $C2^{-k\alpha}$. It
follows that $(2.14)$ is bounded by a constant multiple of
$$ \sum_{k,n} 2^{n+k} \sum_{l} \frac{1}{l} 2^{-n(1-\eta) \alpha}
2^{-k \alpha} \lesssim 1 \tag2.15$$ if $(1-\eta)\alpha>1$, since
$l$ runs up to $C2^{n(1-\eta)-k}$, the number of sectors.

If $j=j'$, $(2.14)$ takes the form
$$ \sum_{k,n} 2^{n+k} \sum_{j'} \int_S
d\mu(x)d\mu(y)d\mu(x')d\mu(y')\tag2.16$$ which is bounded by the
same argument.

\subhead Case 3: $2^n|x-y|>>1$ and $2^n|x'-y'| \lesssim 1$
\endsubhead This case basically vacuous, which can be seen
as follows. We have
$$ \int \int \int \int e^{2 \pi i
t({|x-y|}_{a,\phi}-{|x'-y'|}_{a,\phi})} t^{-\frac{1}{2}}
{|x-y|}_{a, \phi}^{-\frac{1}{2}} \ t d\mu^{*}\beta(2^{-n}t) dt
\psi(a)da d\phi$$
$$=2^n 2^{\frac{n}{2}} \int_{\{2^n|x'-y'| \lesssim 1\}}
\int \int
\widehat{\beta}(2^n({|x-y|}_{a,\phi}-{|x'-y'|}_{a,\phi}))
{|x-y|}_{a, \phi}^{-\frac{1}{2}} \ \psi(a)da d\phi d\mu^{*}$$
$$ \lesssim 2^n 2^{\frac{n}{2}} \int_{\{2^n|x'-y'| \lesssim 1\}}
\int \int \widehat{\beta}(2^n({|x-y|}_{a,\phi})) {|x-y|}_{a,
\phi}^{-\frac{1}{2}} \ \psi(a)da d\phi d\mu^{*}$$ $$ \lesssim 2^n
2^{\frac{n}{2}}2^{-n \alpha} \int \int \widehat{\beta}(2^n(|x-y|))
{|x-y|}^{-\frac{1}{2}} \ d\mu(x)d\mu(y) \tag2.17$$

Localizing to the sets where $2^{-k} \leq |x-y| \leq 2^{-k+1}$, we
obtain
$$ 2^n 2^{\frac{n}{2}}2^{-n \alpha} 2^{\frac{k}{2}}
\int \int \widehat{\beta}(2^n(|x-y|)) \ d\mu(x)d\mu(y) \lesssim
2^n 2^{\frac{n}{2}}2^{-n \alpha} 2^{\frac{k}{2}} 2^{-n \alpha}
\tag2.18$$ which sums since $k<<n$. This completes the proof of
Theorem 0.3, and, consequently, the proof of Theorem 0.1.

\vskip.125in

\head Proof of Theorem 1.2 \endhead

\vskip.125in

Define the measure $\nu_0$ by
$$ \int f d\nu_0=\int f({||x-y||}_{K^{*}}) d\mu(x)d\mu(y).
\tag3.1$$

Let
$$ d\nu(s)=e^{i \frac{\pi}{4}} s^{-\frac{1}{2}}
d\nu_0(s)+e^{-i\frac{\pi}{4}} {|s|}^{-\frac{1}{2}} d\nu_0(-s).
\tag3.2$$

Since $\nu_0$ is supported on $\Delta_{K^{*}}(E)$, $\nu$ is
supported on $\Delta_{K^{*}}(E) \cup -\Delta_{K^{*}}(E)$.

We have
$$ \int_{\partial K} {|\widehat{\mu}(t \omega_K)|}^2 d\omega_K=\int
\widehat{\sigma}_t*\mu d\mu, \tag3.3$$ where $\sigma$ is the
measure on $\partial K$.

Using a variant of Theorem 1.4 (see e.g. \cite{Herz62}), we see
that
$$ \widehat{\sigma_t}(x)=2 {(t\rho^{*}(x))}^{-\frac{1}{2}}
\cos \left(2 \pi \left(t \rho^{*}(x)-\frac{1}{8}\right)
\right)+O({(t|x|)}^{-\frac{3}{2}}), \tag3.4$$ where
$$ \rho^{*}(x)=\sup_{y \in \partial K} x \cdot y. \tag3.5$$

In other words, $\rho^{*}(x)={||x||}_{K^{*}}$.

By definition,
$$ \hat{\nu}(k)=e^{i \frac{\pi}{4}} \int
{||x-y||}_{K^{*}}^{-\frac{1}{2}} e^{-2 \pi i k {||x-y||}_{K^{*}}}
d\mu(x)d\mu(y)$$
$$+e^{-i \frac{\pi}{4}} \int
{||x-y||}_{K^{*}}^{-\frac{1}{2}} e^{2 \pi i k {||x-y||}_{K^{*}}}
d\mu(x)d\mu(y)$$
$$=2 \int {||x-y||}_{K^{*}}^{-\frac{1}{2}} \cos \left(2 \pi
\left(|k| \rho^{*}(x-y)-\frac{1}{8}\right) \right)d\mu(x)d\mu(y).
\tag3.6$$

By $(3.3)$,
$$ \int {|\widehat{\mu}(k \omega_K)|}^2 d\omega_K={|k|}^{-\frac{1}{2}}
\int 2{||x-y||}_{K^{*}}^{-\frac{1}{2}} \cos \left(2 \pi \left(|k|
\rho^{*}(x-y)-\frac{1}{8}\right) \right) d\mu(x)d\mu(y)$$
$$+O \left( \int_{|x-y| \ge {|k|}^{-1}}
{(|k||x-y|)}^{-\frac{3}{2}} d\mu(x)d\mu(y) \right)$$
$$+O \left(\int_{|x-y| \leq {|k|}^{-1}}
{(|k||x-y|)}^{-\frac{1}{2}} d\mu(x)d\mu(y) \right)$$
$$={|k|}^{-\frac{1}{2}}
\int 2{||x-y||}_{K^{*}}^{-\frac{1}{2}} \cos \left(2 \pi \left(|k|
\rho^{*}(x-y)-\frac{1}{8}\right) \right) d\mu(x)d\mu(y)$$
$$+O \left( \int {(|k||x-y|)}^{-\frac{1}{2}} d\mu(x)d\mu(y)
\right) \tag3.7$$ for any $\alpha \in [1/2,3/2]$. It follows that
$$ \hat{\nu}(k)={|k|}^{\frac{1}{2}} \int {|\widehat{\mu}(k
\omega)|}^2 d\omega+O({|k|}^{\frac{1}{2}-\alpha}I_{\alpha}(\mu)).
\tag3.8$$

Since the error term is clearly in $L^2(|k| \ge 1)$, we see that
$\hat{\nu} \in L^2(|k| \ge 1)$ if and only if ${|k|}^{\frac{1}{2}}
\int {|\widehat{\mu}(k \omega_K)|}^2 d\omega_K \in L^2(|k| \ge
1)$. This precisely what Theorem 1.2 asserts.

\newpage

\head References \endhead

\vskip.125in

\ref \key Bourgain94 \by J. Bourgain \paper Hausdorff dimension
and distance sets \jour Israel J. Math. \vol 87 \yr 1994 \pages
193-201 \endref

\ref \key Falconer86 \by K. J. Falconer \paper On the Hausdorff
dimensions of distance sets \jour Mathematika \vol 32 \pages
206-212 \yr 1986 \endref

\ref \key Herz62 \by C. Herz \paper Fourier transforms related to
convex sets \yr 1962 \jour Ann. of Math. \vol 75 \pages 81-92
\endref

\ref \key IoLa2003 \by A. Iosevich and I. Laba \paper $K$-distance
sets and Falconer conjecture \jour (in preparation) \yr 2002
\endref

\ref \key Mattila87 \by P. Mattila \paper Spherical averages of
Fourier transforms of measures with finite energy: dimensions of
intersections and distance sets \jour Mathematika \vol 34 \yr 1987
\pages 207-228 \endref

\ref \key PA95 \by J. Pach and and P. Agarwal \paper Combinatorial
Geometry \yr 1995 \jour Wiley-Interscience Series \endref

\ref \key Sogge93 \by C. D. Sogge \paper Fourier integrals in
classical analysis \jour Cambridge University Press \yr 1993
\endref

\ref \key Sjolin93 \by P. Sjolin \paper Estimates of spherical
averages of Fourier transforms and dimensions of sets \yr 1993
\jour Mathematika \vol 40 \pages 322-330 \endref

\ref \key Wolff02 \by T. Wolff \paper Lectures in Harmonic
Analysis \jour California Institute of Technology Class Lectures
Notes (revised by I. Laba) \yr 2002 \endref

\enddocument